\documentclass[12pt]{article}
\usepackage{amsmath,amsfonts, amssymb,verbatim}
\newtheorem{pkt}{}[section]
\newcommand{\bpk}{\begin{pkt}\rm }
\newcommand{\epk}{\end{pkt}}

\newcommand{\Uu}{\mathrm{U}}
\newcommand{\Vv}{\mathrm{V}}

 \newcommand{\E}{\mathrm{E}}

\newcommand{\U}{\mathbb{U}}
\newcommand{\uu}{{\bf u}}
\newcommand{\vv}{{\bf v}}

\newcommand{\Oo}{\mathrm{O}^\mathrm{Def}}
\newcommand{\hh}{\mathbf{h}}
\newcommand{\M}{{\mathbb{M}}}

\newcommand{\nn}{\mathfrak{n}}

\newcommand{\Pp}{\mathrm{P}}
\newcommand{\Qq}{\mathrm{Q}}

\newcommand{\R}{{\mathbb R}}
\newcommand{\Q}{{\mathbb Q}}
\newcommand{\Z}{{\mathbb Z}}
\newcommand{\N}{{\mathbb N}}
\newcommand{\C}{{\mathbb C}}

\newcommand{\be}{\begin{equation}}
\newcommand{\ee}{\end{equation}}

\newcommand{\HH}{\mathcal{H}}

\newcommand{\invdx}{\sqrt{\frac{\nn}{2\pi}}}
\newcommand{\dx}{\Delta_\nn}

\newcommand{\ra}{\rangle}
\newcommand{\la}{\langle}

\newcommand{\e}{\mathrm{e}}

\begin{document}
\title{Dirac - von Neumann axioms in the setting of Continuous Model Theory }
\author{B.Zilber}

\maketitle
%\tableofcontents

\begin{abstract}
We recast the well-known axiom system of quantum mechanics used by physicists (the Dirac calculus) in the language of Continuous Logic. For the basic version of the axiomatic system we prove that along with the canonical continuous model the axioms have approximate finite models of large sizes, in fact the continuous model is isomorphic to an ultraproduct of finite models. We analyse the continuous logic quantifier corresponding to Dirac integration and show that in finite context it has two versions, local and global, which coincide on Gaussian wave-functions.

\end{abstract}

\section{Introduction}

\bpk 

The axiomatic formulation of quantum mechanics was introduced by Paul Dirac in 1930 \cite{Dirac} through a description of a Hilbert space, and later developed with greater mathematical rigor in a monograph of 1932 by John von Neumann. Since 1930, has gone through several rewritings and new editions, refining the calculus to a level he considered satisfactory.  {In fact, von Neumann himself expressed
	his dissatisfaction not long after publication of his book,
	and spent considerable effort looking for alternatives \cite{Redei}.} 
%Modern books present Dirac's axioms in a succinct form, often omitting much of the technical detail.

This discussion continues presently, see e.g.  \cite{Nonphysical}, which  questions the physicality of the complete Hilbert space $\HH$ of quantum mechanics, effectively arguing for a dense subspace $\HH^\mathrm{Def}$ of physically meaningfull states instead.  
\epk
\bpk 
In Section \ref{s2} we survey   the  axioms of quantum mechanics following \cite{ModernP}.
Readers with a background in logic  will notice that what physicists refer to as axioms is very far from what is a conventional set of axioms  in a formal language even in its early form as presented e.g. by Hilbert's axiomatisation of geometry \cite{HilbertGrundlagen}. 
We noted earlier, in \cite{axiomsCL} and in \cite{ZCL}, that the language that Dirac introduced is, in effect, that of  continuous logic, CL, and the (rigged) Hilbert space which the axioms describe is a CL-analogue of cylindric algebra of Tarski, see \cite{Monk},  \cite{axiomsCL}. 

More precisely, we consider a separable Hilbert space $\HH^{\otimes k} $ as a completion of  a linear subspace $\HH^\mathrm{Def}_k$ generated by $k$-ary continuous predicates $\psi: U^k\to \C$. This correspondence between Hilbert spaces and continuous structures gives us two 
 kinds of models,  which we call $\HH$-structures:

- a canonical model $\U(\infty)$  on $\R,$

-  finite  models  $\U(n)$ forming an asymptotic class

In Section \ref{s4} we prove our main  Theorem \ref{Thm1}, stating that 
the canonical continuous model $\U(\infty)$  is isomorphic to a CL-ultraproduct of  $\U(n),$ $n\in \N.$

An important point about the setting: the CL-model $\U(\infty)$ of Dirac - von Neumann axioms is compact, and so 
its theory is totally (absolutely) categorical, see e.g. \cite{GL}. This places the theory very high in CL-classification hierarchy.  In contrast,  CL-theories of Hilbert spaces with non-commuting operators are unstable; see \cite{FH}.  
 \epk 
\bpk 
It should be noted that, to represent  elements of the rigged Hilbert space which are outside $L^2(\R),$ including all-important position and momentum eigenvectors, one has to use the model-theoretic techniques of  imaginary elements and imaginary sorts, see \ref{genpred}. Moreover, this technique allows us to immediately generalise the setting of $L^2(\R^m)$, of quantum mechanics over the space $\R^m$, to  the general $L^2(\M),$ of quantum mechanics over a manifold $\M$ which is interpretable in $\R.$ 

Note the key difference between the Hilbert space-representation of quantum mechanics and the CL-models representation: the latter has an explicit geometric universe $\M$, similar to the universe of Newtonian physics. 
\epk 
\bpk A consequence of Theorem \ref{Thm1} is that 
 $\U(\infty)$ is pseudo-finite in the sense of continuous model theory. This fact is known in general for  compact structures, see \cite{GL}.
However, our $\U(n)$ are specific lattice models, similar 
to lattice models  in use in physics.  Such constructions play an important role as toy models, as vehicles for particular calculations, and for modelling particular phenomena. However,  
the fact that the specific lattice models  $\U(n),$ hosting the full set of integral operators, represent  Dirac's axioms of quantum mechanics in a mathematically exact way is novel. 
In particular, the Dirac integral is shown to correspond to a summation formula which is ``local'' in a certain well-defined sense, namely the summation domain has to be much shorter then the (non-standard) length of the full universe.

The proof of the theorem is based on an analysis of Dirac's integral (equivalently of the structure of the underlying rigged Hilbert space) in the context of a CL-quantifier. Modelling Dirac's integral in  CL-structures $\U(n)$  exposes a specific  local nature of this quantifier. In fact, the approximation construction summarised by the main theorem \ref{Thm1} is exactly the local approximation for this class of structures as defined in \cite{ZLG}.

 In contrast, it is possible to approximate the canonical continuous model by  another type of approximation,  called   
  global in \cite{ZLG}, provided we restrict the language to Gaussian states (of the form $\e^{Q(x)},$ where $Q$ is a quadratic form). In this approximation the metric can be ignored and the structure is seen globally as a discrete structure.  
 Our second main result, Theorem \ref{Thm2}, states that for the Gaussian theories  the local and the global approximations produce the same continuous structure.

Gaussian part of quantum mechanics is the backbone of the theory. This makes, in our opinion, the noted correspondence between local and global especially interesting.

\epk

\section{Axiomatisation of  quantum mechanics}\label{s2}

Below we reproduce a slighly edited version of a typical representation of axioms of QM for physicists; see \cite{ModernP}, 6.3.

\bpk
{\bf Axiom 1}. The state of a quantum system is described by a vector $|\psi\ra$ of norm 1  belonging to
a complex Hilbert space $\HH.$ This state is usually called ``ket $\psi$''.
A complex Hilbert space $\HH$ is a vector space, which can be finite dimensional
or infinite dimensional, equipped with the complex scalar product (also called inner
product) $\la \psi |\psi' \ra$ between any pair of states $|\psi\ra,$ $|\psi'\ra$ in $\HH.$ The norm, or modulus, of
a generic vector $|\psi\ra\in \HH$ is defined as
$$||\psi||=\sqrt{|\la\psi  \psi\ra}$$
and usually $|\psi\ra$ is normalized to one, i.e.$||\psi||=1.$ The symbol $\la \psi|$ which appears in
the definition of the norm is called ``bra $\psi$'' and it can be intepreted as the function $$\la \psi|: \HH\to \C.$$
 For any $|\psi'\ra\in \HH$  this function gives a complex number $\la \psi|\psi'\ra$ obtained as
scalar product of $|\psi\ra$ and $|\psi'\ra$. In a complex Hilbert space $\HH$ it exists a set of basis
vectors $|\phi_\alpha\ra$  which are orthonormal, 
i.e. $\la \phi_\alpha|\phi_\beta\ra =\delta(\alpha-\beta),$ and such that \be \label{psi}|\psi\ra=\sum_\alpha c_\alpha |\phi_\alpha\ra,\ee
for any $|\psi\ra$, where the coefficients $c_\alpha$ belong to $\C.$

\medskip

{\bf Axiom 2}. Any observable (measurable quantity) of a quantum system is described
by a self-adjoint linear operator $F: \HH\to \HH$ acting on the Hilbert space of state vectors.

For any classical observable $F$ there exists a corresponding quantum observable $F$.

\medskip

{\bf Axioms 3}.The possible measurable values of an observable $F$ are its eigenvalues $f,$
such that
$$ F|f\ra=f |f\ra,$$
with $|f\ra$ the corresponding eigenstate.
The observable $|f\ra$ admits the spectral resolution
\be \label{f}F=\sum_f f |f\ra \la f|,\ee 
where  $\{ |f\ra\}$ is the set of orthonormal eigenstates of $F$, and the mathematical object $\la f|$, called ``bra of $f$'', is a linear map that maps into the complex number. This
also satisfies
 the identity
$$\sum_f |f\ra \la f|=\mathrm{I}.$$

\medskip

{\bf Axiom 4}. The probability $P$ of finding the state $|\psi\ra$ in the state $|f\ra$ (both of norm 1) is given by
$$P = |\la f|\psi\ra|^2.$$
 This probability $P$ is also the probability of measuring the value $f$ of
the observable $F$ when the system is in the quantum state  $|\psi\ra$.

\medskip

{\bf Axiom 5}. The time evolution of states and observables of a quantum system with
Hamiltonian $\mathrm{H}$ is determined by the unitary operator
$$K^t:=\exp(-i\mathrm{H}t/\hbar),$$   
such that $|\psi(t)\ra  = K^t|\psi\ra$ 
is the time-evolved state $|\psi\ra.$ 
\epk
\bpk\label{2.2}  {\bf The dynamical reformulation of quantum mechanics}. This is based on  the Stone Theorem: 

{\em For each self-adjoint operator $A$ on $\HH$ there is a well-defined one parameter group of unitary operators on $\HH$
$$\{ \e^{i At}: \ t\in \R\}$$
and $A$ can be recovered uniquely from the group.}

Thus, we may reduce the  theory to the equivalent theory of Hilbert spaces with unitary operators of the form above. 
One advantage of such a theory is that the unitary operators, unlike unbounded self-adjoint operators, are defined on the whole of $\HH$ and their treatment is  mathematically more straightforward. The framework is also called the Heisenberg picture of quantum mechanics.  
\epk
\bpk  Now we make several comments on the axioms.

The term ``Hilbert space'' here should actually be read  as {\em the rigged Hilbert space} (see \cite{rigged}) because it differs from the standard definition by accommodating both a Hilbert space $\Phi$  and the dual space $\Phi^*$ with $$\Phi\subseteq \HH\subseteq \Phi^*,$$ 
see further details in \ref{rigged}.

The summation formulas like (\ref{psi}) and (\ref{f}) are presented in a form of an integral if the family $|\psi_\alpha\ra$  is continuous but seems natural in the summation form when  $\alpha$ runs in the discrete spectrum of an operator. 

{\em 
Note that our exposition of the basic formulas of Quantum Mechanics in the next section is not meant to follow a well-developed mathematical framework such as the classical Schwartz theory of distributions. Instead, we give an overview of basic QM as practised by physicists, in order to compare it with the presentation of the alternative CL-formalism of $\HH$-structures in Sections 4 and 5. }

\epk
\bpk\label{R3.4} {\bf Remark.} 
Hilbert space formalism provides a powerful mathematical framework to extend quantum mechanics, allowing distributions and generalised eigenfunctions to be rigorously handled. However, as is almost generally accepted, not every element corresponds to a physically realisable state -- some are purely mathematical artefacts; see e.g. \cite{Nonphysical}. 

In the more general context of quantum field theories Wightman axioms explicitly postulate that physically meaningful part of the rigged Hilbert space $\HH$ is a dense subset $\HH^D\subset \HH.$ Our representation takes $\HH^D$ to be exactly the subspace of definable states.
\epk
\section{The axioms in the setting of Continuous Logic}\label{s3}
\bpk\label{CLA} 
Our interpretation of the Axioms  is based, following \cite{axiomsCL}, on the formalism of Continuous Logic (CL), in a slightly generalised form. We do not claim its full equivalence to the Hilbert space formalism but rather demonstrate that standard calculations can be carried out in CL-setting, which may be more natural for physical theories. 

The ideas and the main techniques of continuous model theory is presented e.g. in \cite{Hart}. Its connection to a Hilbert space and Banach space formalisms in the general setting is discussed in \cite{axiomsCL}.

Recall that, in its most general terms, the language of CL consists of {\em predicate} symbols (we will ignore function symbols for now), a collection of {\em connectives}, that is, continuous functions $\C^n\to \C$, and quantifiers, which in general are operators $Q$ on the space of formulas, $\psi(x,\bar y)\mapsto Qx\,\psi(x,\bar y)$.
% For our case, we choose to define these as integration operators; see below and in Section \ref{s4}.

A basic CL-formula is built from predicate symbols using connectives and quantifiers.

An interpretation of symbols and formulas begins with a choice of a {\em universe} $M$, which may be a metric space or, in more recent applications, a metric space equipped with a measure. In the setting of the current paper we assume that $M$ is a union of compact subdomains.  

Symbols of $n$-ary predicates $P$ are interpreted as maps $$P\colon M^n\to D_P\subset \C,$$ where $D_P$ is a compact domain in $\C$.
If $f\colon \C^n\to \C$ is a connective and $\psi_1,\ldots,\psi_n$ are formulas, equivalently, definable predicates, then the formula $f(\psi_1,\ldots,\psi_n)$ is interpreted as the composition of the maps defined by $\psi_1,\ldots,\psi_n$ with $f$.
In the setting of this paper, the only quantifiers we use are integrals over bounded subdomains of $M$; see details in Section \ref{s4}.

The uniformity of the interpretation of language symbols across different $M$ is ensured by uniform continuity moduli for the symbols $P$. These occur automatically on compact subdomains of $M$.

A universe $M$ together with an interpretation of the predicates $P$ of the language constitutes a {\em structure} in continuous model theory. Importantly, the definable sets in a structure are obtained not only via CL-formulas, but also as limits of families of formula-definable sets which makes e.g. the Riemann integral over a closed interval definable using standard CL-quantifiers. 

The techniques of quantum mechanics require the use of linear functionals over the spaces of predicates. We realise such functionals in the form of generalised predicates, which can also be viewed as imaginary elements. Again, this does not necessarily lead to the full dual space of the Gel'fand triple, but it does deliver position and momentum eigenstates and, with these, a dense subspace of the Gel'fand triple.

\epk
\bpk
Recall that the historical prototype of a vector $|\psi\ra$ of the Hilbert space has been a wave-function, that is a function
$$\psi:  \mathcal{M}\to \C$$
from a configuration space $\mathcal{M}$ into a bounded domain of the complex numbers $\C.$ 

These can be seen as  predicates on a domain which, as the matter of fact is identified in Dirac calculus of quantum mechanics with $\R^n,$ where $\R$ is the real line  seen as a measure space. It is clear that in the setting of continuous logic $\R^n$ can be always reduced to $\R,$ which is then our main structure of interest.

 Definable predicates of norm 1 will be referred to as states. 
Of special significance are the {\bf momentum and position states} with the caveat that their norm is defined to be 1 in the sense of the delta-function.  Momentum states were defined by Dirac as the definable family of predicates of the form
\be\label{p} |p\ra:= \frac{1}{\sqrt{2\pi }}\e^{-ipx},\ee
for $p\in \R,$
and they form an orthonormal system in the sense of the inner product

{ 
%??One can consider a $\C$-linear space generated by the momentum states and
%define Hermitian inner product, first between the momentum states
\be\label{pp} \la p_1|p_2\ra:=\delta(p_1-p_2),\ee
where $\delta$ is the Dirac delta.

The  {\bf position} states $|x\ra,$ $x\in \R,$ take the role of characteristic functions of one-point subsets $\{ x\}$ 
\be \label{px1}|x\ra:= \delta(x-z)\ee
(as a function of $z$),
which for  convenience of continuous mathematical manipulations   have been replaced here and in (\ref{pp}) by   Dirac's delta-functions, that is by distributions. In this sense one writes
\be \label{xp} |x\ra=\frac{1}{\sqrt{2\pi}}\int_\R \e^{ixp}|p\ra dp,\ee

Equivalently, position states can be represented by linear functionals (bra-vectors)
$$\la x|:\  |\psi\ra\mapsto \psi(x).$$

Applying the  Fourier duality one can represent $|x\ra$ in the dual configuration space $\{ p\in \R\}$ as

\be\label{xx} |x\ra:= \frac{1}{\sqrt{2\pi }}\e^{ipx},\ee
as functions of $p,$ in analogy to (\ref{p}). Respectively, in analogy to (\ref{px1}) 
$$|p\ra:= \delta(p-t)$$
as a function of $t\in \R.$

%or equivalently, for $\psi$ running in $\{ |p\ra\ : \ p\in \R\}$,
%\be \label{pxF}\la x|:\  |p\ra\mapsto \frac{1}{\sqrt{2\pi }}\e^{-ipx}.\ee
}

We will show that $|p\ra$ and $|x\ra$ can be interpreted in the context of continuous logic as imaginary elements.

The basic unitary operators (Weyl operators) can be defined by their action on the basis:
$$\e^{i\Pp}: |p\ra \mapsto \e^{ip}|p\ra,$$
$$\e^{i\Qq}: |x\ra \mapsto \e^{ix}|x\ra.$$
In particular, the former can be equivalently, using (\ref{xp}), written as
$$\e^{i\Pp}: |x\ra \mapsto \frac{1}{\sqrt{2\pi}}\int_\R \e^{ip(x+1)}|p\ra dp.$$

The time-evolution operator $K^t_\mathrm{Free}$ for a free particle is 
$\e^{-it\frac{\Pp^2}{2}},\ t\in \R,$ that is 
\be \label{p2} K^t: |p\ra \mapsto \e^{-it \frac{p^2}{2}}|p\ra,\ee  
which yields by (\ref{xp})
$$ K^t: |x_0\ra \mapsto \frac{1}{\sqrt{2\pi}}\int_\R \e^{i(px_0-t \frac{p^2}{2})}|p\ra dp$$ 
and
$$\la x| K^t |x_0\ra=\frac{1}{\sqrt{2\pi}}\int_\R \e^{i(px_0-t \frac{p^2}{2})}\la x|p\ra dp.$$ 
Substituting (\ref{p})   one gets
$$ \frac{1}{\sqrt{2\pi}}\int_\R \e^{i(px_0-t \frac{p^2}{2})}\la x|p\ra dp=
\frac{1}{{2\pi}}\int_\R \e^{i(p(x_0-x)-t \frac{p^2}{2})} dp=$$
$$
\frac{1}{{2\pi}}\int_\R \e^{-it \frac{(p-\frac{(x-x_0}{t})^2}{2}+\frac{i(x-x_0)^2}{2t}}dp  = \frac{1}{{2\pi}}\e^{i\frac{(x-x_0)^2}{2t}}\cdot \left(\int_\R \e^{-it \frac{p^2}{2}} dp\right)=$$
$$=\frac{1}{{2\pi}}\e^{i\frac{(x-x_0)^2}{2t}}\cdot \sqrt{\frac{2\pi}{it}}=
\frac{1}{\sqrt{2\pi it}}\e^\frac{i(x-x_0^2)}{2t},$$
and one obtains the well known formula determining the probability amplitude of the free particle propagating from position $|x_0$ to position $|x\ra$ in time $t$:
 \be \label{p2t} \la x|K^t |x_0\ra = \frac{1}{\sqrt{2\pi it}}\e^\frac{i(x-x_0)^2}{2t}.\ee

\epk
\bpk In a more abstract axiomatic setting the theory of Dirac integration in the context of Gaussian states lead to the following formulae:

for $a\neq 0,$  
\be \label{GaussIntegral} \int_\R \e^{\pi i ({a x^2+ 2bx})} dx=\sqrt{\frac{ 1 }{ia}}\e^{-\pi i\frac{b^2}{a}}=\e^{-\frac{\pi i}{4}}\sqrt{\frac{ 1}{a}} \e^{-\pi i\frac{b^2}{a}}, 
\ee
 and, for $a=0,$ 
\be \label{Gpx}  \int_\R \e^{-2\pi i bx} dx= \,\delta(b),\ee
which are not proper Riemann integrals.  

 These formulas will be given an interpretation in  CL-setting.

\epk
\bpk
In terms of structures,
let  $\HH_m$ be the set of all $m$-ary predicates on $\R.$ This by definition has structure of $\C$-vector spaces $$\C=\HH_0\subset \ldots \subset\HH_m\subset\ldots \subset\HH_{m+1}\ldots \HH.$$   

Also, one uses quantifiers, linear maps written as integrals 
$$\phi(z_1,\ldots,z_n)\mapsto \int_\R \phi(z_1,\ldots,z_n) dz_n.$$
In fact, this is a collection of linear maps 
$$\int:\ \HH_{m+1}\to \HH_m,$$ 
the rules of calculation of which  defined by Dirac's  improper integration.

For a discrete basis the sum in (\ref{psi}) can be represented in CL-setting as a definable countable sum, see e.g. \cite{Hart}, Example 5.2.

A special binary operation in the spaces,
 inner product, 
 $$\HH_m\times \HH_m\to \C; \ \ \la \phi(z_1,\ldots,z_n), \psi(z_1,\ldots,z_n)\ra = \int_{\R^m}\phi\cdot \psi^c\, dz_1\ldots dz_m,$$
where $\psi^c$ is the complex comjugate of $\psi$ and $\int_{\R^m}$ is the m-fold integral.
 $\la \phi|\psi\ra$
   can be seen as a continuous predicate of equality $\phi=\psi.$ 

   One restricts the notion of {\bf state} to predicates $\phi$ such that $\la \phi|\phi\ra=1.$  
   
   \medskip

An important role in the theory is played by a collection of linear maps (operators)  
$$L:  \HH_m\to \HH_m$$
with physical meanings.
These can be of the integral  form 
$$\phi(\bar{z}_1,\bar{z_2})\mapsto \int_{\R^k} \alpha(\bar{y},\bar{z_1})\cdot  \phi(\bar{y},\bar{z}_2)\ d \bar{y},$$
where $ |\bar{y}|=|\bar{z}_1|=k,$ $\alpha\in \HH_{2k},$ 
or as classical linear operators
 $$\Pp: \phi(x,\bar{z})\mapsto i\hbar \frac{\partial \phi(x,\bar{z})}{\partial x}\ \mbox{ or } \ \Qq: \phi(x,\bar{z})\mapsto x\cdot\phi(x,\bar{z}).$$ 
  
 The {\bf time evolution operator} $\exp(-i\mathrm{H}t/\hbar)$  acts on states as a unitary operator determining the evolution of a state in time $t$ with a given Hamiltonian H. A state $\phi_{t_0}$ determining a system  at time $t_0$ evolves into a state $\phi_{t}:=\exp(-i\mathrm{H}(t-t_0)/\hbar)$ with the {\em probability amplitude} equal to
 $\la \phi_{t_0}|\phi_{t}\ra,$ which is a complex number of modulus 1.
 The calculation of the CL-formulae $\phi_{t}$ and $\la \phi_{t_0}|\phi_{t}\ra$ (which involve mainly calculations of the application of quantifier $\int$) is the central problem of quantum theory, equivalent to solving the associated Schr\"odinger equation.

 All of the above together makes the $\HH_m$ a collection of Hilbert spaces with  linear operators and  $\HH$ an ambient Hilbert space.
 
 \epk
 
\bpk \label{dynamic}
Below we explain how the Hilbert space axiomatisation of QM can be represented as a formal theory in the language of Continuous Logic. We may assume that operators on the Hilbert space are all integral kernel operators, see \ref{2.2} above, although other forms of operators are also definable in our structures, see \ref{ops}.

\epk

\bpk \label{rigged}
{\bf Remarks on rigged Hilbert spaces.} 

Recall (see e.g. \cite{Bohm})  that  a Gel'fand triple is:
$$\Phi \subset \HH\subset \Phi^*,$$
where $\Phi$ is a space of test functions (e.g. the space of continuous functions on $\R$ with compact support), 
$\HH$ 
is the Hilbert space, $\Phi^*$
is the space of continuous linear functionals on $\Phi$, i.e., Schwartz distributions.

In the context of continuous model theory with universe $\R$, it is natural to take for $\Phi$ continuous functions which are zero outside the interval $[-m,m],$ for some $m,$ in the Banach space topology.
This agrees  well with the conventional approach.   

An element $\phi\in \Phi^*$ acts on a test function $f$  
via the application of inner product in $\Phi^*$: $f\mapsto \la f|\phi\ra$

The ket $|x\ra \in \Phi^*$
is not a vector in 
$\HH$, but a generalised eigenvector.

The pairing $\la x|\phi\ra$
can be interpreted as the evaluation $\phi(x)$ of $\phi$
at the point $x$, if such evaluation makes sense. 
If $\phi\in \HH,$ 
 this is a well-defined function $\phi: \R\to \C.$
%in the sense that, though not necessarily continuous or point-wise defined everywhere.
\epk 
\section{The Hilbert space formalism and $\HH$-structures }\label{s4}
\bpk The axiomatic description of quantum mechanical theory in the form of rigged Hilbert space 
may be quite confusing from the logician point of view -- there are no logical sentences which can be  called axioms.
What Axioms 1 -- 5 render instead is the topological-algebraic structure of a Hilbert space with operators. 
This brings us to the {\em algebraisation of logic} approach
%, perhaps less popular among model theorists nowadays, 
introduced by A. Lindenbaum, A. Tarski and P. Halmos for the first order setting.
It is quite natural to see the Hilbert space formalism as the form of algebraic logic in the context of the continuous logic of physics.

As explained in \cite{axiomsCL}, the Hilbert space $\HH$ of quantum mechanics, or rather the tower of Hilbert spaces $\HH^{\otimes n},$ plays the role of the Tarski {\em cylindric algebra}, see \cite{Monk}.

 \epk

\bpk\label{defHstr0} {\bf  $\HH$-structure.}% (in place of \ref{defHstr}).

Let $\HH$ be a rigged Hilbert space, $\HH_n=\HH^{\otimes n}$ and $\Oo$ a collection of linear operators on the $\HH_n.$ We will often write $\HH$ for the union of the tower $\HH_1\subset \HH_2\subset \ldots $ %and similarly with $\HH^\mathrm{Def}.$

Let $\HH_n^\mathrm{Def}$ be a  subspace of the $\HH_n$ closed under operators from $\Oo.$ 

Given $(\HH^\mathrm{Def},\Oo)$ we will associate with it a structure $(\U; \HH^\mathrm{Def},\Oo)$ in continuous model theory (sometimes called an $\HH$-structure).   $\U$ is the metric universe represented as a union of finite diameter domains,  $(\HH^\mathrm{Def},\Oo)$ is the signature.

There is a natural measure on $\U$ which allows the definition of an integral $\int \psi(x) dx$. 
Using the integral we define   
a Hermitian inner product on  the linear spaces $\HH^\mathrm{Def}_n.$    This gives us a structure of continuous model theory with the integral as a quantifier.

Then we can use the construction of imaginary elements and imaginary sorts to recover in $\U$ elements of spaces dual to $\HH^\mathrm{Def}_n.$
% in particular position and momentum eigenstates and respective continuous bases as imaginary sorts.

We aim to represent in $(\U; \HH^\mathrm{Def},\Oo)$ enough of physically meaningful states, at least all-important position $|x\ra$ and momentum $|p\ra$ eigenvectors, and $\Oo$ to contain  
general integral kernel operators, which includes   Weyl operators $\e^{i\Qq}$ and $\e^{i\Pp},$
the time evolution operators for basic particles and many others.

\epk
\bpk \label{defHstr}
{\bf The construction.} The {\bf universe} of $\U,$ which we also call $\U,$ is a complete metric space  with a distinguished point $0,$ a measure $\mu$ and  metric  $\mathrm {d}(u_1,u_2).$

There are two cases of how the universe can be presented. 
In case $\U=\U(\infty)$ is infinite, define  $\U(\infty):=\R,$ $\mathrm {d}(u_1,u_2)=|u_1-u_2|$ and the measure is the standard measure on $\R.$
%Dirac's delta-measure as determined in \ref{RemarkIntegration}.

In case $\U=\U(\nn)$ is finite of size $\nn$ (may assume $\nn$ even and $\sqrt{\nn}$ is an integer) we identify
$$\U(\nn):=[-\frac{\nn}{2},\frac{\nn}{2})\cap \Z $$
with metric given by   
$\mathrm {d}(u_1,u_2)=\dx\cdot
|u_1-u_2|$, 
 where $|u_1-u_2|$ is the natural integer valued distance on $\Z$ and $\dx=\sqrt{\frac{2\pi}{\nn}}.$
 The measure of a point is $dx=\dx.$

The  universe $\U$ contains  the  family of closed  subdomains $[-m,m]$ of finite diameters, $m\in \Q_{>0},$   %around  the distinguished point $0$
$$ [-m,m]=\{  u\in \U: \mathrm {d}(-m,u)+\mathrm {d}(u,m)=\mathrm{d}(-m,m) \},$$
and we assume that the end-points of the intervals are named.  

In $\U(\infty)=\R$ we consider Lipschitz-continuous, on each subdomain, $n$-ary predicates  
$$\psi: [-m_1,m_1]\times \ldots \times [-m_k,m_k]\to \C.$$ Such predicates are identified with elements of $\HH^{\otimes k}$ and generate the subspace   $\HH^\mathrm{Def}_k$ (we often omit the subscript). There is a natural uniform continuity modulus for each name $\psi.$

Equivalently, we can associate with each such name the map $\psi': \U^k\to \C$ which coincides with $\psi$ on $ [-m_1,m_1]\times \ldots \times [-m_k,m_k]$ and is 0 outside.
Call such $\psi'$ {\bf test functions} with name $\psi.$

The vocabulary of $\HH$-structures $\U$ also includes the global function symbol for a shift $u\mapsto u+\hh.$ The symbol $\hh$ is interpreted in $\U(\nn)$ by an integer $h_\nn$ and $\U(\infty)$ by a real number $h_\infty$ so that
\be \label{+h} \mathrm{d}(u,u+\hh)=\dx \cdot h_\nn \mbox{ in } \U(\nn) \mbox{ and }
 \mathrm{d}(u,u+\hh)=h_\infty \mbox{ in } \U(\infty)  \ee
and the subdomains are shifted by $\hh$: 
\be \label{m+h}  [ -m, m] +\hh= [ -m+\hh,m+\hh ].\ee

For each $n,$ let $\Phi_n$ be  the $\C$-linear  space of functions $\phi: \U^n\to \C$ which are representable as finite linear combinations of test functions.  

In case of the discrete $\U(\nn)$ there is a natural embedding 
\be\label{embed} [-m,m]^{\U(\nn)}\subset [-m,m]^{\U(\infty)}; \ u\mapsto \dx\,u, \ee
and thus we define the predicate with the name $\psi$ on $\U(\nn)$ as the restriction of $\psi$  to $[-m,m]^{\U(\nn)}.$ This is obviously extendable to $n$-ary predicates.

We may use notation  $\Phi_n(\infty),$  $\Phi_n(\nn),$  for spaces of predicates on $\U(\infty)$ and  $\U(\nn)$ respectively.

It is well known, for $\U(\infty)=\R,$ that $\Phi_n$ is a dense inner product subspace of $L^2(\R^n),$ the latter in many cases is taken to be  $\HH$ of the axioms. Moreover, 
$$\Phi_n\subset \HH\subset \Phi^*_n,$$
where $\Phi^*_n$ is the dual to $\Phi_n,$ that is the inner product space of anti-linear functionals on $\Phi_n.$

We introduce  {\bf quantifier} $\E: \Phi_{n+1}\to \Phi_n$ uniformly on all domains $[-m,m]^{n+1}$
$$\E_x: \psi(x,y_1,...,y_n) \mapsto \int_\R   \psi(x,y_1,...,y_n) dx=\int_{-m}^m   \psi(x,y_1,...,y_n) dx.$$

In case of discrete $\U(\nn)$  the quantifier is expressible by the summation formula
$$\E_x: \psi(x,y_1,...,y_n) \mapsto \sum_{-m\le k\dx <m}    \psi(k\dx,y_1,...,y_n)\cdot \dx,$$
where 
$k$ runs in $\Z\cap [-\frac{\nn}{2}, \frac{\nn}{2}).$ 

The quantifier allows the definition of the inner product on $\Phi$ (with the natural extension to $\Phi_n$)
\be \label{innPhi}\la \psi|\phi\ra := \E_x\, \psi(x)\cdot \phi(x)^c, \ee
where $\phi(x)^c$ stands for complex conjugation of $\phi(x).$ 

This implies that $\psi$ can be identified with the anti-linear functional $$\phi\mapsto \la \psi| \phi\ra,$$ an element of $\Phi^*_n.$

\medspace

Below, we omit the subscript in $\Phi_n$ and $\Phi^*_n$ when it does not lead to ambiguity. $\E_x$ in this setting is the quantifier over $x\in \U^n.$  
\epk

\bpk\label{genpred} {\bf Generalised predicates and $\Phi^*$.}
In a slight {\bf extension} of the standard continuous model theory setting we consider the continuous $n$-ary maps into a compact subdomain of $\C,$
 $$\psi: \U^n\to D_\psi\subset \C $$ 

such that

\be \label{psiglobal}\psi=\bigcup_{m\in \N} \psi_{| [-m,m]^n}\ee
the union of basic  continuous predicates.

We call these {\bf  generalised predicates} and consider them definable. 
%  This   passes the main test of definablity. It follows from the representation (\ref{psiglobal}):  

\medspace

Note that: { Given  a  generalised predicate $\psi$ on  $\U$ and  an ultrafilter $\mathcal{D},$ the ultraproduct $\psi^\mathcal{D}$ of $\psi$  in the ultrapower $\U^\mathcal{D},$ is a  generalised predicate.}

\medspace

On the other hand, for a  generalised $\psi$ the inner product $\la \phi|\psi\ra$ with $\phi\in \Phi$ is well-defined since the effective domain of integration is bounded.
So $\psi$ can be identified with the anti-linear functional on $\Phi$
$$ \phi\mapsto \la \phi|\psi\ra,$$ 
that is an imaginary element $|\psi\ra\in \Phi^*.$

We define {\bf momentum eigenfunctions} $\vv[p],$ $p\in \R$ as   generalised continuous predicates on $\U(\infty)=\R$  
$$\vv[p]: x\mapsto \frac{1}{\sqrt{2\pi}}\e^{\frac{ixp}{\hh}}.$$

Another way of representing $\vv[p]$ is by considering the binary predicate $\la\vv|\uu\ra$ on $\U(\infty):$
$$\la\vv|\uu\ra: \R^2\to \C; \ \ (x,p)\mapsto \la x|p\ra= \frac{1}{\sqrt{2\pi}}\e^{-\frac{ixp}{\hh}}$$

This defines the family 
$$\mathbb{V}=\{ \vv[p]: p\in \U  \}$$ as an imaginary sort (in the sense slightly broader than that in \cite{Hart}) with a metric given by 
$$\mathrm {d}(\vv[p_1], \vv[p_2]):= \mathrm {d}(p_1,p_2).$$ 

It is clear that $\mathbb{V}\cong \U$ as metric spaces.  $\mathbb{V}$ is an imaginary sort in $(\U; \HH^\mathrm{Def}).$

Define, for $x\in \U,$ the maps

$$ \uu[x]: \mathbb{V}\to \C; \ \vv[p]\mapsto \vv[p](x)= \frac{1}{\sqrt{2\pi}}\e^{\frac{ixp}{\hh}}. $$ 
These are  generalised continuous predicates on sort $\mathbb{V}.$ Moreover, for each $x$,
$\uu[x]$ defines  the unique anti-linear functional on $\Phi$ 
\be\label{uu1} |x\ra: \Phi \to \C; \ \psi\mapsto \psi(x)=\la \psi| \uu[x]\ra. \ee 

This determines $\uu[x]$ (or $|x\ra$ in bra-ket notations) as imaginary elements and $$\U':= \{ \uu[x]: x\in \U  \},$$ 
a sort of imaginaries in a natural isomorphism with $\U.$ The $\uu[x]$ are not in $\HH$ but in Dirac calculus one defines their norm to be the symbolic Dirac $\delta(0).$ In the linear space $\Phi^*$ one can equivalently introduce the inner product
$$\la \uu[x]| \uu[y]\ra:= \delta_{x,y} \  (\mbox{Kronecker } \delta).$$

In this sense $\U'$ is an orthogonal basis for  $\Phi$  in the sense that, for each $\psi \in \Phi,$ 
\be\label{uu} \psi= \E_x  \psi(x) \uu[x]\ee
since, by linearity of $\E,$ $$ \la \, \E_x  \psi(x) \uu[x]\, | \, \uu[y] \,\ra =\psi(y).$$

\epk

\bpk \label{ops}
{\bf Operators.} There are several definable ways of introducing operators on the space  $\Phi$ of test-functions.

1. {\bf Multiplication operators.} Let $\phi$ be a generalised $n$-ary predicate on $\U.$ Define the multiplication operator on $\Phi_n,$
$$M_\phi:  \psi \mapsto \phi\cdot \psi.$$

\medspace

2. {\bf Integral kernel operators.} Let $\kappa: \U^{2}\to \C$ be a  generalised predicate in variable ${z},x.$ Then for all 
$\psi\in \Phi_n$ in variables $x,\bar{y}$   
$$L_\kappa: \psi\mapsto \E_x\,  \kappa({z},x)\cdot\psi(x,\bar{y}) $$%=\varphi(z,\bar{y})$$
is well-defined and determines a linear operator $\Phi_{n}\to \Phi_{n}.$

\medspace

3. {\bf Operators defined by maps.}  These are operators acting on generalised states assuming there are definable continuous maps $f: \U^n\to \U^n.$ 

For each generalised $n$-ary predicate $\phi$
$$T_f: \phi\mapsto \phi\circ f.$$
It is easy to check that such operators preserve inner product and so are unitary.

%{\bf  The Weyl operator}

\epk
\bpk {\bf Examples}

1. For $f$,  the shift function $u\mapsto u+\hh$,  we get as $T_f$ the Weyl operator 
\be \Vv: \uu[r]\mapsto \uu[r+\hh].\ee

A dual Weyl operator is defined as a multiplication operator:

%\Uu: \uu[r]\mapsto \e^{2\pi i  \frac{r}{\sqrt{\nn}}}\cdot \uu[r]\mbox{ and }
%\Vv: \uu[r]\mapsto \uu[r+h_\nn], \mbox{ if }\U=\U(\nn)\end{array}

%We assume that $\Oo$ contains a pair of   (Weyl) operators $\Uu$ and $\Vv$   on $\Phi$ 
%defined  on the  $\uu[r],\ r\in \U:$ 
\be \label{Uuu} %\begin{array}{ll}
\Uu: \uu[r]\mapsto \e^{ i  r\nu}\cdot \uu[r],\ee 
where 
$\nu=2\pi,$ if $\U=\U(\infty)$, and $\nu=\dx^2=\frac{2\pi}{\nn},$ if $\U=\U(\nn).$

The following commutation relation then holds:
\be\label{UV}\Uu \Vv= \e^{ i \nu\hh} \Vv\Uu.\ee

Using (\ref{uu}) we then get 
$$\begin{array}{ll}\Uu \psi= \E_x \psi(x)\cdot\e^{ i \nu x} \uu[x],\\
\Vv \psi= \E_x \psi(x-\hh) \uu[x].\end{array} 
$$ 
Thus $\Uu$ transforms function $\psi(x)$ on $[k,m]$ to $\psi(x)\cdot \e^{ i\nu  x}$ and
$\Vv$ shifts  $\psi(x)$ on $[k,m]$ to $\psi(x-\hh)$ defined as a continuous test function on $[k+\hh,m].$ (Note that test function $\psi(x)$  vanishes outside $[k,m]$).

On the dual set of vectors $\vv[p],\ p\in \U:$
$$ \Vv: \vv[p]\mapsto \e^{- i\nu  p\hh} \vv[p] \mbox{ and } \Uu: \vv[p]\mapsto \vv[p-1].$$

\medspace

2. An important fact of quantum mechanics is that   unitary operators of the form $\e^{iL},$ where $L$ is self-adjoint
can  be represented in the integral form with respective kernels $\kappa,$  when $L$ is a polynomial of the basic operators $\Pp$ and $\Qq,$ see \cite{Hormander}.

\medspace

3. The time evolution operator $\e^{-iHt/\hbar}$ for the quantum harmonic oscillator of frequency $\omega$ can be represented by an integral operator with kernel
$$\sqrt{\frac{\omega}{2\pi i \hbar \sin\omega t}} \exp i\omega\frac{(x_0^2+x^2)\cos \omega t -2x_0x}{2\hbar\sin \omega t},$$ 
see \cite{ZeidlerII}.

\medspace

4. The multiplication operator $$\Qq: \psi \mapsto \mathrm{d}(0,x)\cdot \psi$$
for $\psi(x)$ on $[0,m]$ can play the role of the  position operator of quantum mechanics. This  does not satisfy the definition of \ref{ops} since $\mathrm{d}(0,x)$ is unbounded on $\U$, and so not a generalised predicate. A similar difficulty with including unbounded operators into the Hilbert space formalism arises in quantum mechanics.

It is possible to avoid the use of unbounded self-adjoint operators by working instead with unitary operators of the form $\e^{iL}$ as in the remark 2 above, see \ref{2.2}.

\epk

\bpk  From now on we will consider $\Oo$ to be a set of linear operators on $\HH^\mathrm{Def}$ which are listed in \ref{ops}.

With this definition an $\HH$ structure $(\U, \HH^\mathrm{Def},\Oo)$ can be reformulated as a structure  $(\U, \HH^\mathrm{Def}).$
\epk
\bpk\label{Prop1} {\bf Proposition.} {\em For the continuous $\HH$-structure $(\U(\infty);\HH^\mathrm{Def})$, 
	there is a unique  Hilbert space $\HH$ which contains a dense inner product subspace of continuous predicates on bounded domains  isomorphic to $\HH^\mathrm{Def}$ with definable operators described in \ref{ops}.
	
	Given the $\HH$-structure  
	$(\U(\nn);\HH^\mathrm{Def})$ constructed in \ref{defHstr}, the space  of  predicates on bounded domains with respect to the inner product is an  $\nn$-dimensional  Hilbert space with definable operators described in \ref{ops}. }

{\bf Proof.} This is the content of \ref{defHstr} - \ref{ops}.$\Box$
\epk
\bpk {\bf Other metric universes.} Suppose $\M\subseteq \R^m=\U(\infty)^m$ is a real manifold definable in  $(\U(\infty);\HH^\mathrm{Def},\Oo)$ along with a metric $\mathrm {d}_M$ and a measure $d^\M x$ on it, compatible with the metric and measure on $\U(\infty).$ Then bounded domains in $\U(\infty)^m$  induce bounded domains  on $\M^n$ and we consider all the definable (in   $(\U(\infty);\HH^\mathrm{Def},\Oo)$)   continuous maps on bounded domains of $\M^n$ to be definable continuous predicates on domains in $\M.$ 

Define $\Phi_n(\M)$ to be the space generated by  all $\psi': \M^n\to \C$ where $\psi'(x)=\psi(x)$ for a definable continuous predicate $\psi$ on a bounded domain $D_\psi\subseteq \M^n,$ when $x\in D_\psi$, and $\psi'(x)=0,$ when $x\notin D_\psi.$  
The measure $d^\M x$ on $\U(\infty)$ induces   the respective   quantifier $\E^\M: \Phi_{n+1}(\M) \to \Phi_n(\M)$. As a result we have a well-defined inner product on  $\Phi_n(\M)$ and the closure of  $\Phi_n(\M)$ in the inner product topology gives us a Hilbert space $\HH_n(\M)$ %isomorphic to $L^2(\M^n)$ by construction
which is contained in the dual space  $\Phi_n(\M)^*.$ This provides us with the rigged Hilbert space of quantum mechanics on the manifold $\M.$

\epk

\bpk {\bf CL-ultraproduct.}  We will work with the asymptotic class of finite-dimensional $\HH$-structures $\U(n)$ of signature $\HH^\mathrm{Def}.$ % described in \ref{defHstr}. 

Note that the predicates of $\HH^\mathrm{Def}$ on $\U(\nn)$  can be written in the form
$$\psi(k,\ldots,k_m)= \psi^{(\infty)}(k_1\dx,\ldots, k_m\dx),$$
where $ \psi^{(\infty)}(x_1,\ldots,x_m)$ is a Lipschitz-continuous function into $\C$ on finite diameter domains in  $\R^m.$  
(The predicate with the same name $\psi$ on the continuous $\U$ is just
$ \psi^{(\infty)}(x_1,\ldots,x_m)$.)

The quantifier $\E$ acting on generaised predicates $\psi$ can be represented as  a family of quantifiers
$$\E:= \{ \E^{(m_1,m_2)}: m_1<m_2\in \Z\},$$ on finite intervals of diameter $m_2-m_1<\invdx$ is defined as follows:%    defined in the two steps:

\be\label{local}\E^{(m_1,m_2)}_k \psi(k,\bar{p})  =\mathrm{st}\left( \dx \sum_{k\dx\ge m_1}^{k\dx\le m_2} \psi(k,\bar{p})\right),\ \mbox{for } \U(\nn), \ee
which corresponds to
\be\label{localR}\E^{(m_1,m_2)}_x \psi(x,\bar{y})  = \int_{m_1}^{m_2}
\psi(x,\bar{y})dx,\ \mbox{for } \U(\infty). \ee

We also consider:

\be \label{local2}	\E^\mathrm{loc}_k \psi(k,\bar{p}):= \lim_{m\to \infty}\E^{(-m,m)}_k \psi(k,\bar{p}),\ee
when the limit exists.

Here, in case $\nn$ is finite, the limit should be understood as the value for the maximal $m$ satisfying $2m\le \invdx.$ 
 
For a pseudo-finite $\nn$ it is assumed that $m$ runs in $\N,$ the standard positive integers. In the context of continuous model theory 	$\E^\mathrm{loc}_k$ is definable in terms of 	$\E^{(-m,m)}$ which we will write as just $\E^{(m)}.$

\medskip

While $\E^{(m)}_k \psi$ is well-defined for all  $\psi,$   	$\E^\mathrm{loc}_k \psi$ might be not, for some $\psi$ for infinite $\nn.$

However, we use notation 	$\E^\mathrm{loc}$ for the family  $\{ \E^{(m)}: m\in \N\}$ when it does not lead to confusion.

\epk 
Recall that the definition of an $\HH$-structure furnishes a definable constant $\hh$, which in the finite-dimensional case can also be expressed in the form $\hh(n)=\frac{h_n}{n},$ see (\ref{Uuu}).    
\bpk\label{Thm1} {\bf Theorem.} {\em 
	Let $\mathcal{D}$ be a non-principal ultrafilter  on $\N$ such that 
	$$\lim_{n/\mathcal{D}} \frac{h_n}{n} =\hh(\infty).$$ 
	Then the continuous model theory ultraproduct $(\U^*,\HH^\mathrm{Def})$ of finite $\HH$-structures $(\U(n),\HH^\mathrm{Def})$  is an $\HH$-structure isomorphic to $(\U(\infty),\HH^\mathrm{Def}).$
	
	Equivalently, for every sentence $\sigma$ and every positive $\epsilon$ there is a subset $D_{\sigma,\epsilon}\in \mathcal{D}$, such that for all $n\in D_{\sigma,\epsilon}$ the value of $\sigma$ on  $(\U(\infty),\HH^\mathrm{Def})$ differs from the value of $\sigma$ on  $(\U(n), \HH^\mathrm{Def})$ by no more than $\epsilon.$
	
}

{\bf Proof.} The metric universe $\U^*$  of the ultraproduct is defined as the union of sorts of finite diameter $2m,$ which are limits along the ultafilter $\mathcal{D}$ of sorts of the same diameter of $\U(n)$. This means that for a limit non-standard number $\nn$ and numbers $k\in \U(\nn)$ we set the  limit point $x=k/_\mathcal{D}$ so that $$\mathrm {d}(0,x)=\mathrm {d}(0,k)/_\mathcal{D}.$$
This brings us to  $$x:=\mathrm{st}(k\dx),$$
(the standard part map). In particular, the interval $[-m\invdx, m\invdx]$ in $\U(\nn)$ corresponds to  the interval $[- m,  m]$ in $\R.$
This also agrees with the definition of predicates $\psi$ on the ultraproduct. 

Next we prove the correspondence for  quantifiers.
It is enough to consider unary $\psi: \R\to \C.$
%First, have to agree that the general unary predicate has the form 
By definition $\psi(k) = \psi^{(\infty)}(k\dx),$ where $ \psi^{(\infty)}(x)$ is  Lipschitz-continuous on each $[-m,m]\subset \R.$ 
%$$\E^{(m)} \psi =\frac{1}{\sqrt{\nn}}\sum_{-m \sqrt{\nn}\le  k< m \sqrt{\nn}}  \psi^{(\infty)}\frac{k}{\sqrt{\nn}})=\frac{1}{\sqrt{\nn}}\sum_j\sum_{-m \sqrt{\nn}\le  k< m \sqrt{\nn}} a_j(\frac{k}{\sqrt{\nn}})^j$$

Claim. Given any positive $\epsilon\in \R,$ 

$$|\,\dx\sum_{-m \le  k\dx < m }  \psi^{(\infty)}(k\dx)\ - \ \int_{-m}^m  \psi^{(\infty)}(x)\, dx\, |<\epsilon.   $$

Indeed, the discrete formula is a Riemann sum with spacing $\Delta x=\dx.$ By the left Riemann sums estimate for an interval $(a,b)$
$$\mathrm{Err}\le M_\psi\frac{(a-b)^2}{2N}=M_\psi\frac{(2m)^2}{4m\sqrt{\nn}}$$
where $N$ is the number of  points between $a=m$ and $b=-m$ and  $M_\psi$ the Lipschitz bound on $[a,b].$ Clearly, $M_\psi\frac{(2m)^2}{4m}\in \R$ and thus $M_\psi\frac{(2m)^2}{4m\sqrt{\nn}}<\epsilon$ because $\dx$ is a non-standard infinitesimal. 

Thus the application of the quantifier in the asymptotic class agrees with the quantifier in the ultraproduct. It follows that the inner product operation  in the asymptotic class agrees with the inner product operation   in the ultraproduct, once it is determined by integration. 

Note that the latter is enough to obtain the correspondence for the construction of interpretable linear functionals and the rigged Hilbert space in the asymptotic class and in the ultraproduct. It follows  that the inner product operation $\la \psi_1|\psi_2\ra$ is preserved by the ultraproduct for all generalised predicates $\psi_1,\psi_2\in \HH^\mathrm{Def}.$
$\Box$
\epk

\section{Gaussian states}
\bpk\label{Gstates} {\bf Gaussian predicates.} 
	Call a generalised $m$-ary predicate $\psi(k_1,\ldots,k_m)$ on $\U(\nn)$ basic Gaussian if there is a $\eta\in \C$ and  a positive-definite  quadratic form $Q(x_1,\ldots,x_m)$ over $\Q$ such that 
	$$ \psi(k_1,\ldots,k_m)= \eta\cdot\e^{-\pi i\frac{Q(k_1,\ldots,k_m)}{\nn}}=\eta\cdot\e^{- i\frac{Q(k_1\dx,\ldots,k_m\dx)}{2}}$$
	
	For the continuous $\U(\infty):$ 
	$$ \psi(x_1,\ldots,x_m)= \eta\cdot\e^{ -  i \frac{Q(x_1,\ldots,x_m)}{2}},$$ where $Q$ is  over   $\R.$ 

Since by definition $\e^{- i rp}$ is a Gaussian predicate, we consider the Fourier-dual one point characteristic function $\uu[r]$  to be a Gaussian state. 

\medskip

Note that  $Q(x_1,\ldots,x_m)$ can be written in the form that singles out a particular variable, say $x_1,$
$$
Q(x_1,\ldots,x_m)=a x^2+2x b(\bar{y})+ c(\bar{y})$$
where $x=x_1,$ $\bar{y}$ is the rest of the variables,   $b(\bar{y})$ a linear form and $c(\bar{y})$ a quadratic form.

Now a Gaussian predicate can be written as
$$\psi(k, \bar{p})=\eta\cdot \e^{-\pi i \frac{c(\bar{p})}{\nn}}\cdot \e^{-\pi i \frac{ak^2+2kb(\bar{p})}{\nn}}$$
in the discrete setting, and 
$$\psi(x, \bar{y})=\eta\cdot \e^{- i \frac{c(\bar{y})}{2}}\cdot  \e^{- i \frac{ax^2+2xb(\bar{y})}{2}}$$
in the continuous seting.

\medskip

For the discrete version, if $a\neq 0,$ we call the rational number $a$ the {\bf period of $\psi$} with respect to variable $k.$  

If $a=0$ and $b(\bar{p})=b\cdot (L_1 p_1+\ldots +L_mp_m)$ with $L_1,\ldots,L_m$ a coprime tuple of integers, then 
 the {\bf period of $\psi$} with respect to variable $k$ is equal to $b.$

\epk
\bpk{\bf $\U(\nn)$ in a non-standard setting.}
In this section we will allow $\nn$ to be a non-standard integer, equivalently 
$$\U(\nn)=\prod \U(n)/\mathcal{D},$$ the first-order ultraproduct, for some non-principal ultrafilter $\mathcal{D}$ on $\N.$ 

More precisely, we choose $\mathcal{D}$ so that the non-standard integer $\nn$ is
 {\bf highly divisible},  that is $\nn$  is divisible by all standard integers.

We still may formally use the metric introduced in \ref{defHstr} in non-standard setting.

\medspace

Let $\nn$ be highly divisible and $|\U|=\nn.$ We say that a subset $X\subset \U^m$ is {\bf d-dense} if $X$ contains a submodule of $\U^m$ of finite index. 

\epk
\bpk\label{GLemma} {\bf Lemma} (Gauss summation). {\em On $\U(\nn),$ for $\nn$ highly divisible:
	
%	assuming $a>0$
\be \frac{1}{\sqrt{\nn}}\sum_{-\frac{\nn}{2a}\le k< \frac{\nn}{2a}}
\e^{-\pi i \frac{ak^2+2kb(\bar{p})}{\nn}}=
\sqrt{ \frac{1}{a}}\cdot \e^{\frac{\pi i}{4}}\cdot \e^{-\pi i\frac{b(\bar{p})^2}{a\nn}}\ee
for all $\bar{p}$ in a d-dense subset of $\U^{m-1},$ and equals 0 outside the d-dense subset. 

In case $a=0$ and $b(\bar{p})= p,$  
\be\label{21}\frac{1}{\sqrt{\nn}}\sum_{-\frac{\nn}{2}\le k< \frac{\nn}{2}}
\e^{\pi i \frac{2kp}{\nn}}= \delta^{(\nn)}(p),\ee
where 
$$\delta^{(\nn)}(p)=\left\lbrace 
\begin{array}{ll} 0\ \ \mbox{ if }p\neq 0,\\  \sqrt{\nn} \mbox{ otherwise}.   \end{array}\right. $$ 

}

{\bf Proof.} %Let $p=b(\bar{p})$ 
Let $a=\frac{A}{D}>0$ where $A,D\in \Z.$ We  choose $D$ so that $D\cdot b(\bar{p})$ is over $\Z.$  
Note that $\frac{\nn}{2a}$ is an integer because $\nn$ is divisible by $A$ by assumptions.

Now let $$X_a=\{ \bar{p}\in \U^{m-1}: A| Db(\bar{p}) \}.$$ This is a dense subset of $\U^{m-1}.$ 
For a $\bar{p}\in X_a$  the function 
$$\e^{-\pi i \frac{ak^2+2kb(\bar{p})}{\nn}}$$ of variable $k$ has period $\frac{\nn}{a}$ and 
$\frac{b(\bar{p})}{a}$ is an integer. Thus  the summands in
	 $$ak^2+2kb(\bar{p})= a(k+\frac{b(\bar{p})}{a})^2 -\frac{b(\bar{p})^2}{a}=an^2-\frac{b(\bar{p})^2}{a}$$
	 are integers and we can write 
$$\sum_{-\frac{\nn}{2a}\le k< \frac{\nn}{2a}}
\e^{-\pi i \frac{ak^2+2kb(\bar{p})}{\nn}}	=
 \e^{i\pi \frac{\frac{b(\bar{p})^2}{a}}{\nn}} 
 	\sum_{0\le n<\frac{\nn}{a}}\e^{-\pi i \frac{an^2}{\nn}}= \e^{i\pi {\frac{b(\bar{p})^2}{a\nn}}}\cdot \sqrt{\frac{\nn}{a}} \e^{\frac{\pi i}{4}},
$$
where at the last step we used the classical Gauss' quadratic sums equality. 

In case $\bar{p}\notin X_a$ the Gauss sum is equal $0.$

\medskip

Now consider the case $a=0,$ %$b=\frac{B}{D},$ for $B,D\in \N$ coprime. %Asssume $$\bar{p}\in X_b=\{ \bar{p}\in \U^m(\nn): \ b(\bar{p})\in \Z \}$$ which is obviously a d-dense subset. 
If  $p=0\mod \nn$ then all the summands in (\ref{21}) are equal to 1 and we get  $\sqrt{\nn}$ for the value of the formula. Alternatively, if  $p\neq 0\mod \nn$, then we get all the roots of 1 of order $\nn$   as summands, and so the sum is equal $0.$  
$\Box$
\epk
\bpk {\bf Lemma.} {\em For any d-dense subset $X\subset \U(\nn)^m$ for any $\bar{y}\in \R^m$ there is $\bar{p}\in X$ such that $\mathrm{st}(\dx\bar{p})=\bar{y}.$}

{\bf Proof.} When $p$ runs in $\U(\nn)=[-\frac{\nn}{2}, \frac{\nn}{2}],$ the numbers  $\mathrm{st}(\dx p)$
run continuously between $-\infty$ and $+\infty.$ Thus there is $\bar{p}\in \U(\nn)^m$ such that $\mathrm{st}(\dx\bar{p}).$

Density of $X$ implies that there is a tuple of non-negative standard integers $\bar{d}$ such that  $\bar{p}+\bar{d}\in X.$ But $\mathrm{st}(\dx\bar{d})=\bar{0}$ and thus we can assume $\bar{p}\in X.$ $\Box$
\epk
\bpk\label{gl} {\bf Corollary} (compare with (\ref{GaussIntegral} and (\ref{Gpx}))).

{\em Let $(\U^*,\HH^\mathrm{Def})$ be the ultraproduct constructed in \ref{Thm1}, $\nn$ highly divisible, and
	 $\e^{-\pi iax^2+2xb(\bar{y})}$ a Gaussian predicate on $\U^*,$ where $a,b$ are finite  rational (possibly non-standard). Then, if $a>0,$

	 \be\label{gla}\frac{1}{\sqrt{\nn}} \sum_{-\frac{\nn}{2a}<k\le  \frac{\nn}{2a}}  \e^{-\pi i \frac{ak^2+ 2kb}{\nn}} = \int_\R \e^{- \pi i {ax^2+2xb}} dx\ee
	  and
	 
	 \be\label{glb}\frac{1}{\sqrt{\nn}}\sum_{-\frac{\nn}{2}<k\le  \frac{\nn}{2}}  \e^{- \frac{2\pi i k}{\nn}} = \int_\R \e^{-2\pi i x} dx\ee   
 where we dropped $\delta^{(\nn)}$ on the left and Dirac delta on the right of the equality.
   }

\epk

\bpk {\bf Global versions of quantifiers.} 
The    {\bf global}  quantifier is defined for finite and pseudo-finite $\U:$
$$\E^\mathrm{glob}_k \psi(k) : =\mathrm{st}\left(\frac{1}{\sqrt{\nn}} \sum_{-\frac{\nn}{2P}<k\le  \frac{\nn}{2P}} \psi(k)\right)$$
where $P$ is the period of $\psi.$

\epk
%We say that an ultrafilter $\mathcal{D}$ on $\N$ is a {\bf divisibility ultrafilter} if $d\cdot \N$ is in $\mathcal{D}$ for each $d\in \N.$

\bpk \label{Thm2} {\bf Theorem.} {\em Given an  $\HH$-structure  $\U=\U(\nn),$ with $\nn$ highly divisible, and a Gaussian  predicate $\psi$ with variables $k, \bar{p}$:  
	$$\E^\mathrm{glob}_k \psi(k,\bar{p})= \E^\mathrm{loc}_k \psi(k,\bar{p})$$
	for each $\bar{p}$ in a d-dense subset.
}

{\bf Proof.} This is a direct consequence of \ref{gl} and \ref{Thm1}. $\Box$
\epk

\end{document}